\documentclass{amsart}
\pdfoutput=1
\usepackage{caption}
\usepackage{subcaption}
\usepackage{amssymb}
\usepackage{graphicx}
\usepackage{wrapfig}
\begin{document}

\author{G. Fejes T\'oth}
%\thanks{Supported by Hungarian OTKA grant K76154 and by the Fields
%Institute for Research in Mathematical Sciences.}
\address{Alfr\'ed R\'enyi Institute of Mathematics,
Re\'altanoda u. 13-15., H-1053, Budapest, Hungary}
\email{gfejes@renyi.hu}

\title{Packing and covering in higher dimensions}
\thanks{The English translation of the book ``Lagerungen in der Ebene,
auf der Kugel und im Raum" by L\'aszl\'o Fejes T\'oth will be
published by Springer in the book series Grundlehren der
mathematischen Wissenschaften under the title
``Lagerungen---Arrangements in the Plane, on the Sphere and
in Space". Besides detailed notes to the original text the
English edition contains eight self-contained new chapters
surveying topics related to the subject of the book but not
contained in it. This is a preprint of one of the new chapters.}

\begin{abstract}
The present work surveys problems in $n$-dimensional space with $n$ large.
Early development in the study of packing and covering in high dimensions
was motivated by the geometry of numbers. Subsequent results, such as the
discovery of the Leech lattice and the linear programming bound, which
culminated in the recent solution of the sphere packing problem
in dimensions 8 and 24, were influenced by coding theory. After mentioning
the known results concerning existence of economical packings and coverings
we discuss the different methods yielding upper bounds for the density
of packing congruent balls. We summarize the few results on upper bounds
for the packing density of general convex bodies. The paper closes with
some remarks on the structure of optimal arrangements.
\end{abstract}

\maketitle
%\section{Higher dimensions}

\section{Existence of economic packings and coverings}

The celebrated Min\-kowski-Hlawka theorem ({\sc Hlawka} \cite{Hlawka}) gives the
lower bound
$$
\delta_L(K)\ge\frac{\zeta(n)}{2^{n-1}}
$$
for the packing density of an arbitrary $n$-dimensional centrally symmetric convex body.
For sufficiently large $n$, {\sc Schmidt} \cite{Schmidt63a,Schmidt63b} improved
the bound to
$cn2^{-n}$, provided $c<\ln\sqrt2$. By the observation of Minkowski
concerning translates of an arbitrary convex body $K$ and $\frac{1}{2}(K-K)$
and the inequality
$$
\frac{{\rm{vol}}(K-K)}{{\rm{vol}}(K)}\le\binom{2n}{n}
$$
of {\sc Rogers} and {\sc Shephard} \cite{RogersShephard}, it follows that, for
sufficiently large $n$,
$$
\delta_L(K)\ge \sqrt{\pi}cn^{3/2}/4^{-n},
$$
where $c<\ln\sqrt2$.

Slightly better lower bounds were proved for the packing density of the ball
by {\sc Ball} \cite{Ball92}, {\sc{Vance}} \cite{Vance}, and {\sc{Venkatesh}}
\cite{Venkatesh}. Venkatesh proved that for any constant
$c>\sinh^2(\pi e)/\pi^2e^3=65963.8\ldots$ there is a number $n(c)$ such that
for $n>n(c)$ we have $\delta_L(B^n)\ge cn2^{-n}$. Moreover, there are
infinitely many dimensions $n$ for which $\delta_L(B^n)\ge n\ln\ln{n}2^{-n-1}$.

{\sc{Leech}} \cite{Leech64,Leech67} constructed lattice ball packings
in dimensions $n\le24$ and $n=2^m$. He found a remarkable packing in $E^{24}$.
This lattice, called after its inventor the {\it{Leech lattice}}, has
density $\pi^{12}/12!=0.001929\ldots$, and each ball in it is touched by
196560 others. As it was shown (see later), no other packing in $E^{24}$
has higher density and no ball can be touched by more than 196560 other
balls of the same size. All the densest known lattice packings in dimensions
less than 24 can be obtained as sections of the Leech lattice.

Leech did not mention, possibly he was not even aware of the fact, that all his
packings were extensions of error-correcting codes, that is packings where
the centers are vertices of the unit cube. A systematic study of constructions
of ball packings based on error-correcting codes was initiated by {\sc{Leech}}
and {\sc{Sloane}} \cite{LeechSloane}. Their work inspired further research, which
led to more constructive lower bounds for the packing density of special
classes of convex bodies, approaching, and in some cases also improving the
Minkowski-Hlawka bound.

Let $B_p^n$ denote the $n$-dimensional unit {\it{superball}}, that is, the ball
for the $l_p$-norm $\left({\sum_{i=1}^n}|x_i|^p\right)^{1/p}$.
{\sc{Rush}} and {\sc{Sloane}} \cite{RushSloane} adopted the constructions
by Leech and Sloane for packings of superballs. They obtained an improvement
of the Minkowski-Hlawka bound for all integers $p>2$. As examples let us
mention that $\delta_L(B_3^n)\ge2^{-0.8226n+o(n)}$ and
$\delta_L(B_4^n)\ge2^{-0.6742n+o(n)}$. With further elaboration of the
method {\sc{Rush}} \cite{Rush89} constructed lattice packings with density
$2^{-n+o(n)}$ for every convex body which is symmetric through each of the
coordinate hyperplanes. Moreover, {\sc{Elkies, Odlyzko}} and {\sc{Rush}}
\cite{ElkiesOdlyzkoRush} were also able to use the method for packings of centrally
symmetric convex bodies. This made it possible to construct dense packings of
{\it{generalized superballs}} defined as
$$f(x_1,\ldots,x_k)^\sigma+f(x_{k+1},\ldots,x_{2k})^\sigma+
\ldots,+f(x_{n-k+1},\ldots,x_n)^\sigma\le1\qquad\sigma\ge1,$$
where $f(x_1,\ldots,x_k)$ is a distance function, that is $f(0)=0$, $f(x)>0$
for $x\ne0$, and $f(tx)=tf(x)$ for $t>0$. With a further generalization of
the notion of superball {\sc{Rush}} \cite{Rush93} could also handle bodies
like $(|x_1|+|x_2|^2)^3$ and $\max(|x_1|,|x_2|,|x_3|^5)$. More results
based on the construction of packings through error correcting codes can be found
in {\sc{Rush}} \cite{Rush91,Rush94,Rush96} and {\sc{Liu}} and {\sc{Xing}}
\cite{LiuXing}.

For packing of balls in spherical space, independently of each other,
{\sc{Chabauty}} \cite{Chabauty}, {\sc{Shannon}} \cite{Shannon}, and
{\sc{Wyner}} \cite{Wyner} proved the lower bound
$$
M(n,\varphi)\ge(1+o(1))\sqrt{2\pi{n}}\frac{\cos\varphi}{\sin^{n-1}\varphi},
$$
which was improved by {\sc{Jenssen}}, {\sc{Joos}} and {\sc{Perkins}}
\cite{JenssenJoosPerkins} by a linear factor.

Concerning coverings, {\sc Rogers} \cite{Rogers57} gave the bound
$$
\vartheta_T(K)\le n\ln n +n\ln\ln n+5n\,,
$$
for the translational packing density of an arbitrary $n$-dimensional
convex body $K$. Rogers' proof uses a
mean value argument combined with saturation. There are several
alternative arguments yielding an upper bound for $\vartheta_T(K)$
of the same order of magnitude. {\sc{Nasz\'odi}} \cite{Naszodi16a}
gave a proof of the inequality $\vartheta_T(K)\le n\ln n +n\ln\ln n+5n$
that relies on an algorithmic result of {\sc{Lov\'asz}} \cite{Lovasz}.
{\sc Erd\H{o}s} and {\sc Rogers} \cite{ErdosRogers} showed that every
$n$-dimensional convex body $K$ admits a covering by translates of $K$ with
density $n\ln n +n\ln\ln n+4n$ so that no point is covered more than
$e(n\ln n +n\ln\ln n+4n)$ times. {\sc{F\"uredi}} and {\sc{Kang}}
\cite{FurediKang} used the {\it{Lov\'asz local lemma}} to give a simple proof
of a slightly weaker result: There exists a covering by translates of $K$ such that
every point is covered at most $10{n\ln{n}}$ times.  {\sc{G. Fejes T\'oth}}
\cite{FTG09} observed that with a modification of a proof of {\sc{Rogers}}
\cite{Rogers59} one can show for every $n$-dimensional convex body $K$ the
existence of a lattice arrangement of $K$ such that $O(\ln n)$ translates
of this arrangement form a covering of space with density not exceeding
$n \ln n + n\ln\ln n + n + o(n)$. Thus, a low density covering can be
achieved with an arrangement of relatively simple structure. {\sc{Rolfes}} and
{\sc{Vallentin}} \cite{RolfesVallentin} suggest a greedy approach to
constructing coverings of compact metric spaces by metric balls. Balls are
iteratively chosen to cover the maximum measure of yet uncovered space.
Their method is an extension of the argument of {\sc{Chv\'atal}}
\cite{Chvatal79} for the finite set cover problem to the setting of
compact metric spaces.

{\sc Rogers} \cite{Rogers59} proved that
$$\vartheta_L(K)\le{n^{\log_2\ln n+c}}$$
for some suitable constant $c$ and all $n$-dimensional convex bodies $K$.
The same bound was proved by {\sc{Butler}} \cite{Butler} in a different
way, providing a covering by translates of $K$ with the above density such
that the corresponding translates of $\lambda{K}$,
$\lambda={\rm{vol}}\,{K}/{\rm{vol}}{(K-K)}$, form a packing. For the covering
density of balls {\sc Rogers} \cite{Rogers59} proved the bound
$$\vartheta_L(B^n)\le cn(\ln{n})^{\frac{1}{2}\log_22\pi e}.$$
This result was generalized by {\sc{Gritzmann}} \cite{Gritzmann}, who
showed that
$$\vartheta_L(K)\le{cn(\ln{n})^{1+\log_2e}}$$
for every $n$-dimensional convex body $K$ that has an affine image symmetric
through at least $\log_2\ln n+4$ coordinate hyperplanes.

All the above mentioned results about lattice coverings have been superseded
by {\sc{Ordentlich, Regev}} and {\sc{Weiss}} \cite{OrdentlichRegevWeiss},
who proved that there is a constant $c$ so that for any $n$-dimensional
convex body $K$
$$\vartheta_L(K)\le cn^2.$$

Using covering codes instead of error correcting codes, {\sc{Rush}}
\cite{Rush92} adapted a construction by Leech and Sloane for the construction
of thin lattice coverings by star-shaped bodies.

{\sc Rogers} \cite{Rogers63}, {\sc B\"{o}r\"{o}czky} and {\sc
Wintsche} \cite{BoroczkyJrWintsche03}, {\sc
Verger-Gaugry} \cite{Verger-Gaugry}, and {\sc Dumer} \cite{Dumer07} showed the
existence of reasonably thin coverings of the $n$-dimen\-sional sphere by
congruent spherical caps and an $n$-dimensional ball by congruent balls.
Dumer's result implies $\vartheta(B^n)\le(\frac{1}{2}+o(1))n\log n$ as
$n\to\infty$, improving the above bound of Rogers in the case of a ball
by a factor 1/2. We note that the published version of Dumer's paper
contained a small error which he corrected in \cite{Dumer18}. {\sc{Nasz\'odi}}
\cite{Naszodi16a} established the existence of thin coverings of the
$n$-dimensional sphere by congruent copies of a spherical convex body.

{\sc{Bacher}} \cite{Bacher} considered the problem of finding universal
covering sets.  Let ${\mathcal{F}}$ be a family of convex bodies all having the
same area. A set $U$ is a {\it{universal covering set}} for ${\mathcal{F}}$
if the translates of any member of ${\mathcal{F}}$ by the vectors of $U$
cover the space. He showed the existence of such sets in all dimensions. In
$n$-dimensional space he constructed a universal covering set $U$ for convex bodies of
unit volume such that $rB^n$ contains at most $\ln(r)^{n-1}r^nV(B^n)$ points of $U$.

\section{Upper bounds for $\delta(B^n)$ and lower bounds for  $\vartheta(B^n)$}
\subsection{Blichfeldt's bound}
The first upper bound for the packing density of the $n$-dimensional ball was
proved by {\sc{Blichfeldt}} \cite{Blichfeldt29}. Blichfeldt's idea was the following.
Consider a packing of unit balls in $E^n$. Replace each ball $S$ by a concentric
ball of radius $\sqrt2$ such that the density at distance $d$ from the center of
the ball is $2-d^2$. Of course, the enlarged balls may overlap, however, it can
be shown that their total density is at most 2 at every point of space.
Comparing the volume of the unit ball to the mass of an enlarged ball, we get the bound
$$\delta(B^n)\le\frac{n+2}{2}2^{-n/2}$$
for the packing density of the $n$-dimensional ball. With a modified density
function Blichfeldt obtained a slightly better bound, and a further improvement
was given by {\sc{Rankin}} \cite{Rankin47}. Blichfeldt's method was used by
{\sc{Rankin}} \cite{Rankin55} and {\sc{Bloh}} \cite{Bloh56} to derive bounds for
ball packings in $n$-dimensional spherical space.

\subsection{The simplex bound}
In an $n$-dimensional space of constant curvature let $d_n(r)$ be the density
of $n+1$ mutually touching balls of radius $r$ with respect to the simplex
spanned by the centers of the balls. {\sc{Coxeter}} \cite{Coxeter64b} and
{\sc{L.~Fejes T\'oth}} \cite{FTL53b} conjectured that in an $n$-dimensional
space of constant curvature the density of a packing of balls of radius $r$
cannot exceed the {\it{simplex bound}} $d_n(r)$. The corresponding simplex
bound for coverings, formulated by {\sc{L.~Fejes T\'oth}} \cite{FTL59g, FTL60d},
reads as follows: In an $n$-dimensional space of constant curvature the density
of a covering by balls of radius $r$ cannot be less than the density $D_n(r)$
of $n+1$ balls of radius $r$ centered at the vertices of a regular simplex of
circumradius $r$ relative to the simplex. In Euclidean space the bounds do
not depend on $r$, so in this case we write them simply as $d_n$ and $D_n$.
Concerning the notion of density and the interpretation of the simplex bound
in hyperbolic space see Section 11.2.

For the Euclidean space, the conjecture about packings was verified by
{\sc{Rogers}} \cite{Rogers58} and independently by {\sc{Baranovski\v{i}}}
\cite{Baranovskii64}, and the conjecture about coverings was proved by
{\sc{Coxeter, Few}} and {\sc{Rogers}} \cite{CoxeterFewRogers}.

It can be seen that the simplex bound for $\delta(B^n)$ is sharper than
Blichfeldt's bound. However, the improvement is only a constant factor,
namely we have
$$\delta(B^n)\le d_n\approx\frac{n}{e}2^{-n/2},\quad n\to\infty.$$
The simplex bound
$$\vartheta(B^n)\ge D_n\approx\frac{n}{e\sqrt{e}},\quad n\to\infty,$$
for the covering density of the ball compares quite favorably with the bound
$\vartheta(B^n)\le(\frac{1}{2}+o(1))n\log n$.

The three-dimensional case of the conjecture about packing was settled by
B\"o\-r\"ocz\-ky in the paper by {\sc{B\"o\-r\"ocz\-ky}} and {\sc{Florian}}
\cite{BoroczkyFlorian} and, finally, {\sc{B\"or\"oczky}} \cite{Boroczky78}
proved the conjecture in full generality. More precisely, B\"or\"oczky
proved that the density of each ball in its Dirichlet cell is at most $d_n(r)$.
The conjecture about covering for spherical and hyperbolic space is
still open.

Consider the balls inscribed in the cells of the $3$-dimensional spherical tiling
$\{2,3,3\}$, $\{3,3,3\}$, $\{4,3,3\}$ or $\{5,3,3\}$. The corresponding radii are
${\frac12}\arccos(-\frac13)$, ${\frac12}\arccos(-\frac14)$, $\pi/4$, and $\pi/10$,
and the density of the balls in the Dirichlet cells is $0.61613\ldots$,
$0.68057\ldots$,  $0.72676\ldots$, and  $0.77412\ldots$, respectively. These
densities agree with the corresponding tetrahedral density bound. {\sc{B\"or\"oczky,
B\"or\"oczky Jr, Glazyrin}} and {\sc{Kov\'acs}} \cite{BoroczkyBoroczkyGlazyrinKovacs}
proved the stability of the simplex bound in the cases mentioned here.

Recall that the ordinary sphere cannot be packed
as densely, nor can be covered as thinly by at least two congruent circles as
the Euclidean plane. Remarkably, in three dimension, the analogous statement
does not hold for packings, and probably does not hold for coverings either.
Namely, the density
$$
\frac{60}{\pi}\left(\frac{\pi}{5}-\sin\frac{\pi}{5}\right)=0.77412\ldots\,
$$
of the 120 balls inscribed in the cells of the tiling $\{5,3,3\}$
is greater than $\pi/\sqrt{18}=\delta(B^3)$.
Similarly, the $120$ balls circumscribed about the same cells form a covering
of the spherical space, with density
$$
\frac{60}{\pi}\left(\omega -\sin\omega \right)=1.44480\ldots\,,\
\omega=\frac{2\pi}{3}-\arccos\frac{1}{4}\,,
$$
which is smaller than $5\sqrt5\pi/24=1.46350\ldots\,$, the conjectured
value of $\vartheta(B^3)$. It is conjectured (see {\sc{L. Fejes T\'oth}}
\cite{FTL59g}) that the inspheres and the circumspheres of $\{5,3,3\}$
form the densest packing and the thinnest covering, respectively with
at least 4 spheres in $S^3$. {\sc Florian} \cite{Florian05}) proved that
$d_3(r)$ is a strictly decreasing function for $0<r\le\arctan\sqrt 2$.
It follows that the density of a packing of at least 4 equal spheres in
$S^3$ is at most
$\lim\limits_{r\to0}\sqrt{18}(\arctan\frac{1}{3}-\frac{\pi}{3})=0.77963\ldots$.
This bound is rather close to the conjectured minimum density.

{\sc Moln\'ar} \cite{Molnar63} defined a {\it Segre-Mahler polytope} as a convex
$n$-dimensional polytope in a space of constant curvature, every dihedral angle
of which is at most $120^\circ$. He conjectured that, when equal spheres
of radius $r$ are packed in such a region of $n$-space, the density cannot
exceed the simplex bound $d_n(r)$. He verified the conjecture for the
3-dimensional case.

{\sc{G.~Fejes T\'{o}th}} \cite{FTG80} observed that if a packing of
$n$-dimensional balls arises as the intersection of a higher-dimensional packing
of congruent balls with an $n$-dimensional subspace, then the simplex bound
$d_n$ for the density of the packing still holds. Obviously, no similar
conclusion can be drawn for an arbitrary covering with balls. However,
{\sc{A.~Bezdek}} \cite{BezdekA84} proved that if a circle covering of the
plane arises from a planar section of a lattice covering with balls in three
dimensions, then the covering's density cannot be smaller than
$2\pi/\sqrt{27}+0.017\ldots$. Moreover, equality occurs for exactly one
lattice and only for certain cutting planes.

It should be mentioned that in general there is no connection between the
density of an arrangement of bodies and the densities of the sections of
the arrangement. {\sc{Groemer}} \cite{Groemer66a} gave examples of packings
${\mathcal{P}}_1$ and ${\mathcal{P}}_2$ with density $0$ and $1$, respectively, such
that in each plane parallel to a given plane the density of the intersection
with the sets of ${\mathcal{P}}_1$ is $1$, while in each of these planes the
density of the intersection with the sets of ${\mathcal{P}}_2$ is $0$. He gave
similar examples for coverings.

{\sc{Florian}} \cite{Florian06} gave a nice survey on the simplex bound
for packings of balls in spaces of constant curvature with emphasis on
dimension 2 and 3.

\subsection{The linear programming bound}

It took more than 40 years until an improvement in exponential order
was achieved for Blichfeldt's bound of $\delta(B^n)$. The first step
was made by {\sc{Sidel\'nikov}} \cite{Sidelnikov73,Sidelnikov74}, who
proved $\delta(B^n)\le2^{-0.509619n+o(n)}$. Subsequently,
{\sc{Leven\u{s}te\u{\i}n}} \cite{Levenstein75} improved the bound to
$\delta(B^n)\le2^{-0.5237n+o(n)}$, and {\sc Kabatjanski{\u\i}} and
{\sc Leven\v{s}te\u{i}n} \cite{KabatjanskiiLevenstein} proved
$$
\delta(B^n)\le 2^{-(0.599+o(1))n}\quad ({\rm{as}}\ n\to\infty),
$$
which remains the best known asymptotic upper bound for $\delta(B^n)$.
The gap between this bound and the lower bound by {\sc{Venkatesh}}
\cite{Venkatesh} remains considerable. It is even unknown whether
$\displaystyle\lim_{n\to\infty}\frac{\ln\delta_L(B^n)}{n}$ and
$\displaystyle\lim_{n\to\infty}\frac{\ln\delta(B^n)}{n}$ exist, and if so,
whether they are equal.

All these improvements of the Blichfeldt bound were obtained as a corollary
of a bound for ball packings in spherical space. A {\emph{spherical
code}} in dimension $n$ with minimum angle $\varphi$ is a set of
points on the unit sphere in $E^n$ with given minimum angular distance
$\varphi$ among them. Let $M(n,\varphi)$ denote the greatest size of such
a spherical code. Equivalently, $M(n,\varphi)$ is the maximum number of
balls of angular diameter $\varphi$ that can be packed on the sphere.
{\sc{Delsarte, Goethals}} and {\sc{Seidel}} \cite{DelsatreGoethalsSeidel}
and independently {\sc Kabatjanski{\u\i}} and {\sc Leven\v{s}te\u{i}n}
\cite{KabatjanskiiLevenstein} adopted the linear programming bound
developed by {\sc{Delsarte}} \cite{Delsarte72} for bounding the
cardinality of error-correcting codes of given minimum distance to
spherical codes.

In the linear programming bound a basic role is played by a sequence of
polynomials $P^n_k$, $k=0,\,1,\ldots$, called {\it{ultraspherical
polynomials}}. They form a complete orthogonal system on the interval
$[-1, 1]$ with respect to the measure $(1-t^2)^{(n-3)/2}dt$. In other
words, for $i\neq j$,
$$\int_{-1}^1P^n_i(t)P^n_j(t)(1-t^2)^{(n-3)/2}\,dt=0.$$
For the purpose of the linear programming bound, the normalization is
irrelevant; the sign should be chosen so that $P^n_k(1)>0$. With different
normalizations, they are special Jacobi polynomials or Gegenbauer
polynomials. The property of ultraspherical polynomials that is
crucial for the application is that they are {\it{positive-definite
kernels}}, that is for any $N$ and any points $x_1,\ldots,x_N\in S^{n-1}$,
the $N\times N$ matrix $\left(P^n_k(x_i\cdot x_j)\right)_{1\le i,j\le N}$
is positive semidefinite. {\sc{Schoenberg}} \cite{Schoenberg} proved
that every continuous positive-definite kernel on $S^{n-1}$ has an
ultraspherical expansion with nonnegative coefficients which converges
absolutely and uniformly. Since the sum of the entries of a
positive-semidefinite matrix is nonnegative, it follows that
$$\sum_{1\le i,j\le N}P^n_k(x_i\cdot x_j)\ge0$$
for any points $x_1,\ldots,x_N\in S^{n-1}$.

Now suppose that $x_1,\ldots,x_N$ is a spherical code in dimension $n$ with
minimum angle $\varphi$ and let $P(t)=\sum_{k=0}^ma_kP^n_k(t)$ be a real
polynomial. Then
$$\sum_{1\le i,j\le N}P(x_i\cdot x_j)=NP(1)+\sum_{1\le i,j\le N,\,i\neq j}P(x_i\cdot x_j)=
N^2a_0+\sum_{k=1}^ma_k\sum_{1\le i,j\le N}P^n_k(x_i\cdot x_j).$$
Suppose that $P(t)\le0$ for $-1\le t\le\cos\varphi$, $a_0>0$, and $a_k\ge0$ for
$i= 1, 2,\ldots,k$. Then ${\displaystyle\sum_{1\le i,j\le N,\,i\neq j}}P(x_i\cdot x_j)\le0$
and $\displaystyle\sum_{k=1}^ma_k\sum_{1\le i,j\le N}P^n_k(x_i\cdot x_j)\ge0$. Hence
we get the linear programming bound for $M(n,\varphi)$:
\smallskip

Let $P(t)={\displaystyle\sum_{k=0}^m}a_kP^n_k(t)$ be a real polynomial such that $P(t)\le0$
for $-1\le t\le\cos\varphi$, $a_0>0$, and $a_k\ge0$ for $i= 1, 2,\ldots,m$. Then
$$M(n,\varphi)\le P(1)/a_0.$$
\smallskip

Let $t_{1,k}^n$ denote the largest root of
$P_k^n(t)$. With the choice of appropriate polynomials
{\sc Kabatjanski{\u\i}} and {\sc Leven\v{s}te\u{i}n}
\cite{KabatjanskiiLevenstein} proved that if
$\cos\varphi\le t_{1,k}^n$ then
$$
M(n,\varphi)\le4\binom{k+n-2}{k}\left(1-t_{1,k+1}^n\right)^{-1}.
$$
From this they derived the asymptotic bound
$$
M(n,\varphi)\le(1-\cos\varphi)^{-n/2}2^{-n(0.099+o(1))}\quad ({\rm as}\
n\to\infty)
$$ for all $\varphi\le\varphi^*=62^\circ\ldots\,$.

{\sc{Sardari}} and {\sc{Zargar}} \cite{SardariZargar} sharpened the
above bound to
$$
M(n,\varphi)\le\binom{k+n-2}{k}\left(1+\frac{2k}{n-1}+\frac{1}{1-t_{1,k+1}^n}\right).
$$
The improvement is by a factor of at most 4 and it does not affect
the exponent in the asymptotic bound.

The asymptotic upper bound for $M(n,\varphi)$ yields a similar bound
for $\delta(B^n)$. Kabatjanski{\u\i} and Leven\v{s}te\u{i}n used
the inequality $\delta(B^n)\le\sin^n(\varphi/2)M(n+1,\varphi)$.
However, there is a slightly better inequality using Blichfeldt's
method. For, it is easily seen that a ball of radius $r=1/\sin(\varphi/2)\le2$
can contain at most $M(n,\varphi)$ centers of a packing of unit balls.
Thus the concentric balls of radius $\sin(\varphi/2)$ form an $M(n,\varphi)$-fold
packing, yielding $\delta(B^n)\le\sin^{n}(\varphi/2)M(n,\varphi)$.
This argument was well known in the community approaching problems
of packing by geometric methods. It was first published by {\sc Leven\v{s}te\u{i}n}
\cite[p. 108]{Levenstein83}. It was rediscovered by {\sc{Cohn}} and {\sc{Zhao}}
\cite{CohnZhao}, who also used it to derive an asymptotic improvement of
the simplex bound in hyperbolic space.

The linear programming bound for spherical codes was used by
{\sc{Od\-lyz\-ko}} and {\sc{Sloane}} \cite{OdlyzkoSloane} and
{\sc{Leven\v{s}te\u{i}n}} \cite{Levenstein79} to solve the
problem of densest packing of $N$ balls in $S^n$ for some
special values of $n$ and $N$. The arrangements
that were characterized as optimal spherical codes in this way
are listed, with the exception of the set of $120$ vertices of
the 600 cell $\{3,3,5\}$, in {\sc{Leven\v{s}hte\u{\i}n}}
\cite[p.~72]{Levenstein92} and \cite[p.~621]{Levenstein98}. The
optimality of the vertices of $\{3,3,5\}$ follows from the
simplex bound, and was proved using the linear programming
bound by {\sc{Andreev}} \cite{Andreev99}.

Elaborations of the linear programming bound by {\sc{Musin}} \cite{Musin06b},
{\sc{Boyvalenkov}} \cite{Boyvalenkov}, {\sc{Pfender}} \cite{Pfender}, {\sc{Bachoc}}
and {\sc{Vallentin}} \cite{BachocVallentin08} and {\sc{Mittelman}} and {\sc{Vallentin}}
\cite{MittelmanVallentin} yielded improved upper bounds for $M(n,\varphi)$ in some low
dimensions. The surveys by {\sc{Boyvalenkov, Dodunekov}} and {\sc{Musin}}
\cite{BoyvalenkovDodunekovMusin}, {\sc{Cohn}} \cite{Cohn10,Cohn17a}, and
{\sc{Viazovska}} \cite{Viazovska18} describe these methods in detail.

{\sc Cohn} and {\sc Elkies} \cite{CohnElkies} and {\sc Cohn} \cite{Cohn02}
modified the linear programming method, obtaining bounds for $\delta(B^n)$
directly. Although they did not improve on the asymptotic bound given above,
their method proved to be exceptionally efficient for $n=8$ and $n=24$.
The method enabled {\sc{Cohn}} and {\sc{Kumar}} \cite{CohnKumar04,CohnKumar09}
to show that the Leech lattice is the unique densest lattice in 24 dimensions
and to give an alternative proof of the result of {\sc{Blichfeldt}} \cite{Blichfeldt35}
about the densest packing of balls in 8 dimensions. Furthermore they proved that no
packing of congruent balls in 24 dimensions can exceed the Leech lattice’s
density by a factor of more than $1+1.65\times 10^{-30}$. A breakthrough was
achieved by {\sc{Viazovska}} \cite{Viazovska17} who succeeded in proving
that $\delta(B^8)=\delta_L(B^8)=\pi^4/384$. Some days later {\sc{Cohn, Kumar, Miller,
Radchenko}} and {\sc{Viazovska}} \cite{CohnKumarMillerViazovska} showed that
$\delta(B)^{24}=\delta_L(B^{24})\pi^{12}/12!$. The paper by {\sc{Cohn}} \cite{Cohn17b}
explains the main ideas leading to this landmark achievement. The
survey and interview by {\sc{de Laat}} and {\sc{Vallentin}}
\cite{deLaatVallentin} is also interesting reading on this subject.

\subsection{Arrangements of points with minimum potential energy}

Given a decreasing potential function $f$ defined on $(0,2]$ and an integer
$N>1$, we wish to place $N$ distinct points $\{x_1,x_2,\ldots,x_N\}$ on the unit
sphere in $n$-dimensional space so that the {\it{potential energy}}
$\sum\limits_{i\neq j}f(|x_i-x_j|)$ is as small as possible. {\sc{Yudin}}
\cite{Yudin} extended the linear programming bound to obtain
an lowáer bound for the potential energy. His result contains the linear
programming bound for $M(n,\varphi)$ as a corollary. Using Yudin's result the
minimum potential energy of 240 points on $S^7$, of 196560 points on $S^{23}$,
and of 120 points on $S^3$ with the potential of $f(x)=x^{n-1}$ on $S^n$ was
determined by {\sc{Kolushov}} and {\sc{Yudin}} \cite{KolushovYudin94} and
{\sc{Andreev}} \cite{Andreev97,Andreev00}. The corresponding optimal
arrangements are the minimal vectors in the $E_8$ root lattice, the minimal
vectors in the Leech lattice and the vertices of the 600-cell $\{3,3,5\}$.

{\sc{Cohn}} and {\sc{Kumar}} \cite{CohnKumar07} succeeded in proving optimality
of certain arrangements not just for a specific, but for a whole class of potential
functions. A real function $f(x)$ is {\it{completely monotonic}} if it is
decreasing, infinitely-differentiable and satisfy the inequalities
$(-1)^lf^{(l)}(x)\ge0$. Typical examples are the inverse power functions
$f(x)=1/x^s$ with $s>0$ and exponential functions $f(x)=e^{-cx}$ with $c>0$.
In Cohn and Kumar's investigation the argument of potential functions is
the square of the distance, rather that the distance itself. They define an
arrangement of points as {\it{universally optimal}} if it is optimal under
every potential function $f(|x-y|^2)$ where $f(x^2)$ is completely monotonic.
Note, that $f(x^2)$ is completely monotonic on an interval $(a,b)$, $a>0$,
then $f(x)$ is completely monotonic on $(a^2,b^2)$, but not vice versa.

{\sc{Cohn}} and {\sc{Kumar}} \cite{CohnKumar07} proved universal optimality of
the sets of vertices of all regular simplicial polytopes in every dimension,
as well as of several other arrangements in dimensions $2$ - $8$ and $21$ -
$24$. These arrangements are listed in Table~1 of their paper. The list
coincides with the list of optimal spherical codes mentioned above. For
$f(x)=x^{-t}$ with $t$ large, the energy is asymptotically determined by the
minimal distance, thus minimizing energy requires maximizing the minimal
distance. Hence, universal optimality of an arrangement implies that it is an
optimal spherical code. With the special choice of the potential function
$f(x)=2-x^{1/2}$ and $f(x)=\ln(4/x)$, respectively, it also follows that these
arrangements maximize the sum, as well as the product, of the distances
between pairs of points. Earlier, {\sc{Andreev}} \cite{Andreev96,Andreev97}
proved that the vertices of the regular icosahedron maximize the product of
the distances between 12 points on $S^2$ and the minimal vectors of the Leech
lattice maximize the sum of the distances between 196560 points on
$S^{23}$. Also, {\sc{Kolushov}} and {\sc{Yudin}} \cite{KolushovYudin97} proved
that the vertices of the regular simplex, the regular cross-polytope and the
minimal vectors of the $E_8$ root lattice maximize both the sum and the product of
the distances  between the corresponding number of points lying on the
respective sphere.

On $S^2$, the only universal optima are a single point, two antipodal points,
an equilateral triangle on the equator, and the vertices of a regular tetrahedron,
octahedron, or icosahedron. That this list is complete follows from a result of
{\sc{Leech}} \cite{Leech57} who enumerated all those configurations on $S^2$ that
are in a (stable or unstable) equilibrium under any force that depends on the
distance only. He showed that the only configurations of this kind are the
vertices of the tilings $\{p,q\}\ \ \left( p,q\ge2, \ \frac{1}{p}+\frac{1}{q}>\frac{1}{2}\right)$
or the face-centers of a tiling $\{2,q\}\ \ (q\ge2)$.

Our knowledge in higher dimensions is limited. It appears that universally optimal
arrangements of points are rare. Through a computer search, {\sc{Ballinger,
Blekherman, Cohn, Giansiracusa, Kelly}}  and {\sc{Schürmann}}
\cite{Ballinger+} found two arrangements of points, one consisting of 40 points
on $S^9$, and another of 64 points on $S^{13}$, which they conjecture to be
universally optimal.

Besides the arrangements of points proved by Cohn and Kumar to be universal
optimal, exact solutions are only known for a few special cases.
{\sc{Dragnev}}, {\sc{Legg}} and {\sc{Townsend}} \cite{DragnevLeggTownsend}
investigated the problem of minimizing the energy under the logarithmic
potential $\log(1/x)$. This is the same problem as finding the maximum of
the product of the distances between $k$ points on $S^n$. They solved the
problem for $k=5$ points on $S^2$. The optimum is attained by the opposite
poles of the sphere and three points distributed evenly on the equator,
that is, by the vertices of a bipyramid. {\sc{Dragnev}} \cite{Dragnev16}
investigated the problem in higher dimensions. He made the conjectured that
the product of the distances between $n+3$ points on $S^n$ attains its maximum
for an arrangement of points that is the union of two mutually orthogonal
regular simplices, one of dimension $\lfloor\frac{n+1}{2}\rfloor$, the other
of dimension $\lfloor\frac{n+2}{2}\rfloor$. He proved the conjecture for
$n=3$ and $n=4$. {\sc{Dragnev}} and {\sc{Musin}} \cite{DragnevMusin} verified
the conjecture for all $n$ by enumerating all stationary configurations of
$n+2$ points on $S^n$ for the logarithmic potential.

Due to applications in stereochemistry, the potential ${\rm{sign}}\,(s)/x^s$
received special attention. Through a computer search, {\sc{Melnyk, Knop}} and {\sc{Smith}}
\cite{MelnykKnopSmith} found conjecturally optimal solutions for up to 16 points
for different values of $s$. For the case of 5 points they conjectured that
there exists a phase transition constant $s_0=15.04808$ so that for $s<s_0$
the triangular bipyramid and for $s>s_0$ a quadrilateral pyramid is optimal.
The conjecture was confirmed in three special cases. The triangular
bipyramid was proved to be optimal by {\sc{Hou}} and {\sc{Shao}}
\cite{HouShao} if $s=-1$ and by {\sc{Schwartz}} \cite{Schwartz13} if
$s=1$ or $s=2$. Subsequently, {\sc{Schwartz}} \cite{Schwartz16,Schwartz20}
determined the exact value of $s_0=15.0480773927797\ldots$, proved
that the triangular bipyramid is optimal for $-2<s<s_0$ and the
quadrilateral pyramid is optimal for $s_0<s\le15+25/512=15.048828125.$
For $s=s_0$ both the triangular bipyramid and the quadrilateral pyramid
minimize the energy. For $s>s_0$ the problem is unsolved.

A further characterization of the triangular bipyramid is due to {\sc{Tumanov}}
\cite{Tumanov} who proved that this arrangement of points on the sphere
constitutes the unique minimizer position under a potential of the form
$f(r)=ar^4-br^2+c$, $a>0$, $b>8a$.

Remarkably, the solution of the sphere packing problem in 8 and 24
dimensions was not the end of the story: {\sc{Cohn, Kumar, Miller,
Radchenko}} and {\sc{Viazovska}} \cite{CohnKumarMillerRadchenkoViazovska}
proved that the $E_8$ root lattice and the Leech lattice are universally
optimal among point arrangements in Euclidean spaces of dimensions 8 and
24, respectively. In other words, they have minimum energy among all
arrangements of points with given density for every potential
function that is a completely monotonic function of the squared distance.

\subsection{Lattice arrangements of balls}

Besides dimensions 3, 8, and 24, where the maximum density of general
packings of congruent balls is known, the value of $\delta_L(B^n)$
was determined for $n=4$ and $n=5$ by {\sc Korkine} and {\sc Zolotareff}
\cite{KorkineZolotareff72,KorkineZolotareff77} and for $n=6$ and $7$ by
{\sc Blichfeldt} \cite{Blichfeldt35}. The lattice covering
density of the ball is known only in dimensions 4 and 5. The case $n=4$
has been established by {\sc{Delone}} and {\sc Ry\v{s}kov} \cite{DeloneRyskov}
and the case $n=5$ by {\sc{Ry\v{s}kov}} and {\sc{Baranovski\u{\i}}}
\cite{RyskovBaranovskii75,RyskovBaranovskii76}.

{\sc{Sch\"{u}rmann}} and {\sc{Vallentin}} \cite{SchurmannVallentin}
designed an algorithm for approximating the values of the lattice covering
density of the ball along with the critical lattice, with arbitrary
accuracy. Implementing the algorithm in dimensions 6, 7, and 8 enabled
them to find the best known lattices for thin sphere covering in these
dimensions.

An {\it${m}$-periodic arrangement} is the union of $m$ translates of
a lattice arrangement. {\sc{Andreanov}} and {\sc{Kallus}}
\cite{AndreanovKallus20} presented an algorithm to enumerate all locally
optimal 2-periodic sphere packings in any dimension, provided there are
finitely many. They implemented the algorithm in 3, 4, and 5 dimensions and
showed that no 2-periodic packing of balls surpasses the density of the
optimal lattices in these dimensions.

\section{Bounds for the packing and covering density of convex bodies}

Finding the packing density $\delta(K)$ for a given convex body $K$, even
finding a meaningful upper bound for it, is generally a very difficult task.
Let $d(K)$ denote the density of the insphere of an $n$-dimensional body $K$
in $K$. Then, obviously, $\delta(K)\le\frac{\delta(B^n)}{d(K)}$, which yields
a non-trivial upper bound for $\delta(K)$ if $d(K)\ge\delta(B^n)$. Using this
bound {\sc{Torquato}} and {\sc{Jiao}} \cite{TorquatoJiao09a,TorquatoJiao09b}
gave non-trivial upper bounds for the packing density of the icosahedron and dodecahedron, as well
as for several Archimedean solids and superballs. The bound obtained in this way
for the octahedron is greater than 1, however we get a nontrivial bound for the
packing density of the regular cross-polytope in dimensions greater than 23.
Moreover, we get that the packing density of the $n$-dimensional regular
cross-polytope approaches zero exponentially as $n$ tends to infinity.

The method of Blichfeldt was used to obtain upper bounds for the translational
packing density of the superball
$B^n_p=\left\{(x_1,\ldots,x_n)\in E^n|\left({\sum_{i=1}^n}|x_i|^p\right)^{1/p}\le1\right\}$
by {\sc{van der Corput}} and {\sc{Schaake}} \cite{vanderCorputSchaake},
{\sc{Hua}} \cite{Hua} and {\sc{Rankin}} \cite{Rankin49a, Rankin49b} (see also
{\sc{Zong}} \cite[Section 6.3]{Zong99b}). For $p\ge2$ their bound was recently
improved by {\sc{Sah, Sawhney, Stoner}} and {\sc{Zhao}} \cite{SahSawhneyStonerZhao}
based on the Kabatjanski{\u\i}-Leven\v{s}te\u{i}n bound for spherical codes.

With an extension of Blichfeldt's method, {\sc G.~Fejes T\'{o}th} and {\sc W.~Kuperberg}
\cite{FTGKuperberg93a} was able to give non-trivial upper bounds for the packing
density of other, suitable convex bodies, {\it e.g.}, for ``longish'' bodies such as
sufficiently long cylinders $B^{n-1}\times[0,h]$ and ``sausage-like'' solids
$B^n+[0,h]$ in ${R}^n$. Applying this method {\sc G.~Fejes T\'{o}th, Fodor} and
{\sc{V\'{\i}gh}} \cite{FTGFodorVigh} gave a non-trivial upper bound for the packing
density of the regular cross-polytope in all dimensions greater than 6.

Elaborating on the method used by {\sc A.~Bezdek} and {\sc{W.~Kuperberg}}
\cite{BezdekAKuperberg90}, {\sc{W. Kusner}} \cite{Kusner14}
obtained another bound for the packing density of finite-length circular
cylinders, namely $\delta(B^2\times[0,h])\le\pi/\sqrt{12}+10/h$. Although
this bound is meaningful (smaller than 1) only if $h>100$ and it improves upon
the bound given in G.~Fejes T\'{o}th and W.~Kuperberg only for $h$ greater
than about $250$, the advantage of Kusner's bound is that it approaches the
packing density of the circle as $h\to\infty$. In a further paper
{\sc W. Kusner} \cite{Kusner16} extended the result of A.~Bezdek and
W.~Kuperberg by showing that the packing density of the set
$B^2\times{R^n}$ is also $\pi/\sqrt{12}$.

We can get upper bounds for the translational packing density of convex bodies
by the observation of Minkowski that a family of translates of a convex body
$K$ forms a packing if and only if the corresponding translates of the
centrally symmetric body $\frac{1}{2}(K-K)$ form a packing. {\sc{Groemer}}
\cite{Groemer61d} proved in this way that $\delta_T(C)\le\frac{2^{n-1}}{2^n-1}$
for every convex cone $C$. For the $n$ dimensional simplex $S_n$ this argument
yields $\delta_T(S_n)\le2^n\binom{2n}{n}^{-n}$. In particular, we have
$\delta_T(T)\le0.4$ for a tetrahedron $T$.

{\sc{Zong}} \cite{Zong14b} proposed a method based on the shadow regions
introduced by {\sc{L.~Fejes T\'oth}} \cite{FTL83} to give upper bounds for the
translative packing density of three-dimensional convex bodies. Applying the
method for the tetrahedron $T$, he established the bound
$\delta_T(T)\le{36\sqrt{10}}/({95\sqrt{10}-4})=0.3840610\ldots$.

{\sc{de Oliveira Filho}} and {\sc{Vellentin}} \cite{deOliveiraFilhoVallentin}
extended the linear programming method to estimate the packing density of
congruent copies of a convex body. {\sc{Dostert, Guzm\'an, de Oliveira Filho}}
and {\sc{Vallentin}} \cite{DostertGuzmandeOliveiraFilhoVallentin} exploited
this method to obtain upper bounds for the translative packing density of
some three-dimensional convex bodies with tetrahedral symmetry, such
superballs and Platonic and Archimedean solids. They improved
Zong's upper bound for the translative packing density of the tetrahedron to
$0.3745$.

In some cases it can be proved that for a body $K$ we have $\delta(K)<1$
or $\vartheta(K)>1$ without establishing a concrete bound. Hlawka conjectured
that the packing density of circular tori cannot be 1. Motivated by this
conjecture {\sc{Schmidt}} \cite{Schmidt61} proved a general theorem which
implies as a corollary, besides the positive answer to Hlawka's conjecture,
$\delta(K)<1$ and $\vartheta(K)>1$ for every smooth convex body.
Schmidt's theorem does not apply for packings of cones.
W. Kuperberg proved that a packing consisting of translates of a cone $C$
and its images $-C$ in $E^3$ cannot have density 1. {\sc{B\'ar\'any}} and
{\sc{Matou\v{s}ek}} \cite{BaranyMatousek} succeeded in proving an explicit
bound smaller than 1 for the density of such a packing.

\section{The structure of optimal arrangements}

In higher dimensions the occurrence of less organized arrangements among the
optimal ones seems to be more frequent. It is likely that the equality
$\delta_L(K)=\delta_T(K)$ fails in dimensions greater than 2, although no
convex body is known for which $\delta_L(K)<\delta_T(K)$. {\sc{Rogers}}
\cite[page~15]{Rogers64} conjectured that $\delta_L(B^n)<\delta_T(B^n)$ for all
sufficiently large $n$. {\sc{Best}} \cite{Best} constructed non-lattice ball packings
in dimensions 10, 11, and 13 that are denser than the densest known lattice
packings. There is a special class of convex bodies for which the equality
$\delta_L=\delta_T$ holds: {\sc{Venkov}} \cite{Venkov} and independently
{\sc{McMullen}} \cite{McMullen} proved that {\it parallelohedra}, that is
those polytopes whose translates tile space also admit a lattice tiling.
We note, that according a theorem of {\sc{Groemer}} \cite{Groemer64b}
parallelohedra are also characterized by the property that space can be
tiled by positive homothetic copies of them.

The equality $\delta(K)=\delta_T(K)$, which holds in the plane for all
centrally symmetric disks, fails already in dimension $3$. There exist
centrally symmetric convex bodies whose congruent copies can pack space
perfectly (tile it without gaps), but whose maximum density attained in a
packing of translates is smaller than 1. One such body is the right
double-pyramid erected over and under the unit square, with height $1/2$.
Moreover, {\sc{A.~Bezdek}} and {\sc{W.~ Kuperberg}}
\cite{BezdekAKuperberg91b} observed that for $n\ge3$, there exist
ellipsoids in $E^n$ for which packing with congruent copies can exceed
the maximum density by translates, that is, the ball packing density.

In their construction A.~Bezdek and W.~ Kuperberg used the theorem of
{\sc{Heppes}} \cite{Heppes60b} that one can place infinite circular cylinders
in the void of every lattice packing of balls. Packing in the cylinders
long ellipsoids of the same volume as the balls we get a mixed packing
of balls and ellipsoids. With a suitable affinity, the balls and
ellipsoids are then transformed into congruent ellipsoids. Refining
this construction, {\sc{Wills}} \cite{Wills} produced a packing of
congruent ellipsoids with density $0.7549\ldots$ and  {\sc{Sch\"urmann}}
\cite{Schurmann02b} constructed dense ellipsoid packings in dimensions
up to 8. The packing of congruent ellipsoids constructed by Sch\"urmann
in $E^8$ exceed the density of the densest packing of balls by more than
$42.9\%$.

Motivated by the problem of understanding the structure of certain materials
like crystals and glasses {\sc Donev}, {\sc{Stillinger}} {\sc{Chaikin}} and
{\sc{Torquato}} \cite{Donev+} constructed in $E^3$ packings of congruent
copies of ellipsoids close to the ball with density grater than
$\pi/\sqrt{18}=0.74048\ldots$. If $a\le{b}\le{c}$ are the semiaxes of the
ellipsoid and $c/a\ge\sqrt3$, then the construction has density 0.7707,
exceeding $\pi/\sqrt{18}$ considerably. In the case $1.365\le{a/b}\le1.5625$
a packing constructed by {\sc{Jin, Jiao, Liu, Yuan}} and {\sc{Li}} \cite{Jin+}
has higher density.

The idea of the combination of a lattice arrangement with an arrangement in
infinite cylinders was used by {\sc{G.~Fejes T\'oth}} and {\sc{W.~Kuperberg}}
\cite{FTGKuperberg95} for coverings. They proved that for every $n$-dimensional
($n\ge3$) strictly convex body $K$ there is an affine-equivalent body $L$ whose
congruent copies can cover space more thinly than any lattice covering. The
assumption of strict convexity is essential: There exist convex polyhedra,
{\it{e.g.}} a cube, which tile space in a lattice-like manner. In $E^3$ they
presented an ellipsoid $E$ for which $\vartheta(E)\le1.394$. Since by
the simplex bound of {\sc{Coxeter, Few}} and {\sc{Rogers}} \cite{CoxeterFewRogers}
$\vartheta(B^3)\ge\frac{3\sqrt3}{2}(3\arccos\frac{1}{3}-\pi)=1.431\ldots$, it follows
that $\vartheta(E)<\vartheta(B^3)=\vartheta_T(E)$.

%Let $p(K)$ and $c(K)$ be the number of, up to isometry, different densest
%lattice packings and thinnest lattice coverings of $K$, respectively. {\sc
%Gruber} \cite{Gruber11} showed that for a typical $n$-dimensional convex body
%$K$ both $p(K)$ and $c(K)$ are bounded by constants depending only on $n$. For
%$n=2$ and $3$ a typical convex body has a unique densest lattice packing and a
%unique thinnest lattice covering.

The excellent book of {\sc{Rogers}} \cite{Rogers64} gives an exhaustive account of
packing and covering in high dimensions.  The book of {\sc{Conway}} and {\sc{Sloane}}
\cite{ConwaySloane} is an encyclopedic source of information about sphere packing.

\small{
\bibliography{pack}}